\documentclass[12pt]{amsart}

\usepackage{latexsym, amssymb}

\newtheorem{thm}{Theorem}
\newtheorem{prop}{Proposition} 
\newtheorem{lma}{Lemma}

\newtheorem{df}{Definition}

\newcommand{\be}{\begin{equation}}
\newcommand{\ee}{\end{equation}}
\newcommand{\beq}{\begin{eqnarray}}
\newcommand{\eeq}{\end{eqnarray}}

\def\S{\Sigma}
\def\Om{\Omega}
\def\SO{\Sigma_O}
\def\SH{\Sigma_H}
\def\tOm{\tilde{\Omega}}

\def\S{\Sigma}

\def\R{\mathbb{R}}

\def\vh{\vspace{.3cm}}

\def\MmR{M^m_o}
\def\Mmrz{M^m_{r_0}}
\def\SR{\Sigma_O}
\def\Sro{\S_\rho}

\def\Sr{\S_r}
\def\bN{\bar{M}}

\def\SMm{ ( M^m, g^m ) }
\def\gm{g_{_M}}
\def\vh{\vspace{.3cm}}
\def\stop{\hfill$\Box$}

\begin{document}

\title{On a Localized Riemannian Penrose Inequality}


\author{Pengzi Miao}


\thanks{This research is supported in part by the Australian Research Council}

\address{School of Mathematical Sciences, Monash University,
Victoria, 3800, Australia.} \email{Pengzi.Miao@sci.monash.edu.au}

\begin{abstract}
Let $ \Omega $ be  a compact, orientable, three dimensional Riemannian 
manifold with boundary  with nonnegative scalar curvature.
 Suppose its boundary $ \partial \Omega $  is the 
disjoint union of two  pieces: $ \SH $  and 
$ \SO $, where $ \SH $ consists of the unique closed minimal 
surfaces in $ \Omega $ and $ \SO $ is metrically a round sphere.
We obtain an inequality relating the area of $ \SH $ to the 
area and the total mean curvature  of $ \SO $.
Such an $ \Omega $ may be thought as a  
region, surrounding the outermost apparent horizons of 
black holes, in a time-symmetric slice of a space-time in the 
context of general relativity. The inequality we establish has 
close ties with the Riemannian Penrose Inequality, proved 
by Huisken and Ilmanen \cite{IMF} and by Bray  \cite{Bray}. 
\end{abstract}

\maketitle

\section{Introduction}
Let $ M $ be a complete, asymptotically flat $ 3 $-manifold 
with nonnegative scalar curvature.  Suppose its boundary 
$ \partial M $ consists of the outermost minimal surfaces in 
$ M $. The Riemannian Penrose Inequality,  first proved by 
Huisken and Ilmanen \cite{IMF} for a connected $ \partial M $, 
and then by Bray \cite{Bray} for  $ \partial M $ with any 
number of components, states that 
\be \label{RPI}
m_{ADM} ( M ) \geq \sqrt{ \frac{ A }{ 16 \pi } } ,
\ee
where $ m_{ADM} ( M ) $ is the ADM mass \cite{ADM} of 
$ M $ and $ A $ is the area of $ \partial M $. 
Furthermore, the equality holds if and only if $ M $ is 
isometric to a spatial Schwarzschild manifold outside its horizon.

Motivated by the quasi-local mass question in general relativity
(see \cite{Bartnik_qlmass}, \cite{BYmass1}, \cite{CDO-Yau}, etc.), 
we would like to seek a localized statement of the above 
inequality (\ref{RPI}). To be precise, we are interested in a 
compact, orientable,  $ 3 $-dimensional Riemannian 
manifold $ \Omega $ with boundary. 
We call $ \Omega $   a {\em body surrounding horizons} 
if its boundary $ \partial \Omega $ is the disjoint union of 
two pieces: $ \Sigma_O $ (the outer boundary) and 
$ \S_H $ (the horizon boundary), and $ \Om $ satisfies 
the following assumptions:
\begin{enumerate}
\item[(a)] $ \S_O $ is topologically a $ 2 $-sphere.
\item[(b)]  Each component of $ \S_H $ is a minimal surface  in $ \Omega $.
\item[(c)] There are no other  closed minimal surfaces in 
$ \Omega $. 
\end{enumerate}
Physically, $ \Omega $ is to be thought as a finite region 
in a time-symmetric slice of  a spacetime containing 
black holes and $ \S_H $ corresponds to the outermost 
apparent horizon of the black holes. In such a context, 
if the spacetime satisfies the dominant energy condition
and $ m_{QM} (\SO) $ represents some quantity (to be defined) which 
could measure the quasi-local mass of $ \SO $, then one would 
 expect  
\be \label{gqlmass}
 m_{QM} (\SO ) \geq  \sqrt{ \frac{ A }{ 16 \pi } } .
\ee

In this paper, we are able to establish an equality of the above form 
for a special class of {\em body surrounding horizons}. Our main result is 

\begin{thm} \label{mthm} 
Let $ \Omega $ be a body surrounding horizons whose 
outer boundary $ \SO $ is metrically a round sphere.
Suppose $ \Om $ has nonnegative scalar curvature and 
$ \SO $ has positive mean curvature. Then
\be \label{mineq}
m ( \SO ) \geq \sqrt{ \frac{ | \SH | }{ 16 \pi } } ,
\ee
where  $ m (\SO) $ is defined by
\be
m ( \SO ) = \sqrt{ \frac{ | \SO | }{ 16 \pi } } 
\left[ 1 - \frac{ 1 }{ 1 6 \pi | \SO | } \left( \oint_{\SO} H \ d \sigma 
\right)^2 \right] ,
\ee
where $ | \SH | $, $ | \SO | $ are the area of 
$ \SH $, $ \SO$, $ H $ is the mean curvature of 
$ \SO $ (with respect to the outward normal) in $ \Om $, 
and $ d \sigma $ is the surface measure of the induced metric. 
When equality holds, $ \SO $ is a surface with constant 
mean curvature. 
\end{thm}

We  remark that, assuming \eqref{mineq}  in Theorem \ref{mthm} holds in the first place, 
one can  derive \eqref{RPI} in the 
Riemannian Penrose Inequality. That is because, by a result 
of Bray  \cite{Bray_Thesis},  to 
prove  \eqref{RPI},  one suffices to prove it for a special
asymptotically flat manifold $ M $ which, outside 
some compact set $ K $, is isometric to a spatial 
Schwarzschild manifold near infinity. 
On such an $ M $,  let $ \Om $ be a  compact region 
containing $ K $ such that its 
outer boundary $ \SO $  is a rotationally symmetric sphere 
in the Schwarzschild region.  Applying Theorem \ref{mthm} 
to such an $ \Om $ and observing that, in this case,  the 
quantity $ m( \SO ) $ coincides with the Hawking quasi-local mass
\cite{Hawking} of $ \SO $, hence agrees with the ADM 
mass of $ M$, we see that  (\ref{mineq}) implies (\ref{RPI}). 
On the other hand, our proof of Theorem \ref{mthm} does make 
critical  use of  (\ref{RPI}).
Therefore, \eqref{mineq} and 
\eqref{RPI} are equivalent. 

Besides the Riemannian Penrose Inequality, 
Theorem \ref{mthm} is also largely inspired by the following 
result of Shi and Tam \cite{Shi-Tam}:

\begin{thm}  \label{stthm} {\em (Shi-Tam) }
Let $ \tOm $ be a compact, $ 3 $-dimensional  Riemannian 
manifold with boundary with 
nonnegative scalar curvature.  Suppose $ \partial  \tOm $ 
has finitely many components $ \S_i $ so that each $ \S_i $ 
has positive Gaussian curvature and positive mean 
curvature $ H $ (with respect to the outward normal), then
\be \label{ST}
\oint_{ \S_i } H \ d \sigma \leq \oint_{ \S_i } H_0 \ d \sigma,
\ee
where $ H_0 $ is the mean curvature of $ \S_i $ 
(with respect to the outward normal) when it is 
isometrically imbedded in $ \R^3 $. Furthermore, 
equality holds if and only if $ \partial \tOm $ has only
one component and $ \tOm$ is isometric to a domain 
in $ \R^3 $.
\end{thm}

Let $ \tOm $ be given  in Theorem \ref{stthm}. 
Suppose $ \partial \tOm $ has a component $ \S $ which 
is isometric to a round sphere with area $ 4 \pi R^2 $, then
\be \label{specialhz}
\oint_\S H_0 \ d \sigma = 8 \pi R ,
\ee
and (\ref{ST})  yields
\be
\frac{ 1 }{ 8 \pi } \oint_{ \S } H \ d \sigma \leq   R .
\ee
Now  suppose there is a closed minimal 
surfaces $ \S_h $ in $ \tilde{\Omega} $ such that $ \S_h $ and $ \S $ 
bounds a region  $ \Om $ which contains no other closed minimal
surfaces  in  $ \tilde{\Om} $ (by minimizing area over surfaces homologous to $ \S$, 
such a  $ \S_h $ always exists if $ \partial \tilde{\Om} $ 
has more than one components).
 Applying Theorem \ref{mthm}
to $ \Om $, we have
\be
\sqrt{ \frac{ | \S_h | }{ 1 6 \pi } } \leq 
\sqrt{ \frac{ | \S | }{ 16 \pi } } \left[ 1 - \frac{ 1 }{ 1 6 \pi | \S | } 
\left( \oint_{\S} H d \sigma \right)^2 \right]  ,
\ee
which can be equivalently written as
\be \label{equiform}
\frac{ 1 }{ 8 \pi } \oint_{ \S } H \ d \sigma \leq \sqrt{ R ( R  - R_h ) },
\ee
where $ R = \sqrt{ \frac{ | \S | }{ 4 \pi }} $ and
$ R_h = \sqrt{ \frac{ | \S_h | }{ 4 \pi } } $.
Therefore, Theorem \ref{mthm}  may be viewed as a 
refinement of Theorem \ref{stthm} in this special case to include the 
effect on $ \S$ by the closed minimal surface in $ \tOm $ 
that lies  ``closest" to $ \S $.

In general relativity, Theorem \ref{stthm} is a statement on the positivity of the Brown-York quasi-local mass $m_{BY}(\partial \tilde{\Om})$ \cite{BYmass1}.  Using the technique of weak inverse mean curvature flow developed by Huisken and Ilmanen \cite{IMF}, Shi and Tam  \cite{Shi-Tam07} further proved that $ m_{BY} ( \partial \tilde{\Om} )$ is bounded from below by the Hawking quasi-local mass $ m_H ( \partial \tilde{\Om} ) $.  Suggested by the quantity $ m ( \SO) $
in Theorem \ref{mthm}, we find some new geometric quantities associated to $ \partial \tilde{\Om}$, which are interestingly 
between $ m_{BY}(\partial \tilde{\Om} ) $ and 
$ m_{H}(\partial \tilde{\Om} )$ (hence providing another  proof of  
 $ m_{BY} (\partial \tilde{\Om})  \geq m_H( \partial \tilde{\Om}  )$.)
We include this discussion at the end of the paper.

This paper is organized as follows. In Section \ref{SST}, we 
review the approach of Shi and Tam in  \cite{Shi-Tam} since it 
plays a key role in our derivation of Theorem \ref{mthm}. The detailed 
proof of Theorem \ref{mthm} is given in Section \ref{mpf}.
In Section \ref{gmono}, we establish a partially generalized
Shi-Tam monotonicity. In Section \ref{apofrpi}, we make use of the 
Riemannian Penrose Inequality. 
In Section \ref{rmk}, we give some discussion on quasi-local mass. In particular, we introduce two quantities motivated by Theorem
\ref{mthm} and compare them with the Brown-York quasi-local 
mass $ m_{BY} ( \Sigma )$ and the Hawking quasi-local mass 
$ m_H ( \Sigma ) $.

\vspace{.3cm}

{\it Acknowledgement}: The author wants to thank professor Hubert Bray for the helpful discussion leading to Section 4.


\section{Review of Shi-Tam's approach } \label{SST}
In \cite{Shi-Tam}, Shi and Tam pioneered the idea of 
using results on asymptotically flat manifolds to study compact
manifolds with boundary. We briefly review their approach in this section.  

 Let $ \tOm $ be given 
in Theorem \ref{stthm}. For simplicity, we assume  
$ \partial \tOm $  has only one component $ \S $. 
Since $ \S $ has positive Gaussian curvature, $ \S $ can 
be isometrically imbedded in $ \R^3 $ as a strictly 
convex surface \cite{Nirenberg}. On the region 
$ E $ exterior to $ \S $, the Euclidean metric $ g_0 $ 
can be written as 
\be
g_0 = d \rho^2 + g_\rho ,
\ee
where $ g_\rho $ is the induced metric on each 
level set $ \Sro $ of the Euclidean distance function 
$ \rho $ to $ \S $. Motivated by the quasi-spherical 
metric construction of Bartnik \cite{Bartnik}, Shi and 
Tam showed that there exists a positive function $ u $ 
defined on $ E $ such that the warped metric 
\be
g_u = u^2 d \rho^2 + g_\rho 
\ee
 has zero scalar curvature, is asymptotically flat 
 and the mean curvature of $ \S $ in $ ( E, g_u )$ 
 (with respect to the $ \infty $-pointing normal) agrees
 with the mean curvature of $ \S $ in $ \Om $.
Furthermore, as a key ingredient to prove their result,  
they showed that the quantity
\be
\oint_{ \S_\rho } ( H_0 - H_u )  \ d \sigma 
\ee
is monotone non-increasing in $ \rho $, and 
\be
\lim_{ \rho \rightarrow \infty } \oint_{ \S_\rho } 
( H_0 - H_u )  \ d \sigma  = 8 \pi m_{ADM}(g_u),
\ee
where $ H_0 $, $ H_u $ are the mean curvature of 
$ \S_\rho $ with respect to $ g_0 $, $ g_u $, and 
$ m_{ADM} ( g_u ) $ is the ADM mass of $ g_u $. 
Let $ M $ be the Riemannian manifold  
obtained by gluing $ ( E, g^m_u )$ to $ \tOm 
$ along $ \S $. The metric on $ M $ is asymptotically flat, 
has nonnegative scalar curvature away from $ \S $, 
is Lipschitz near $ \S $, and the mean curvatures of 
$ \S $ computed in both sides of $ \S $ in $ M $
(with respect to the $ \infty $-pointing normal) are the same. 
By generalizing Witten's spinor argument \cite{Witten}, 
Shi and Tam proved that the positive mass theorem 
\cite{Schoen-Yau}  \cite{Witten} remains valid on $ M $ 
(see \cite{Miao_PMT} for a non-spinor proof).
Therefore, 
\be
\oint_{ \S} ( H_0 - H_u )  \ d \sigma \geq 
\lim_{ \rho \rightarrow \infty } \oint_{ \S_\rho } 
( H_0 - H_u )  \ d \sigma  = 8 \pi m_{ADM}(g_u) \geq 0 ,
\ee
with $ \oint_\Sigma ( H_0 -  H ) \ d \sigma = 0 $
 if and only if  $ H = H_0 $ and 
$ \tOm $ is isometric to a domain in $ \R^3 $. 

\section{Proof of Theorem \ref{mthm}} \label{mpf}

We are now in a position to prove Theorem \ref{mthm}.
The basic idea is to deform the exterior region
of a rotationally symmetric sphere in a spatial 
Schwarzschild manifold  in a similar way as Shi and Tam 
did on $ \R^3 $, then  attach it to a 
{\em body surrounding horizons} and apply the 
Riemannian Penrose Inequality to the gluing manifold. 
The key ingredient in our proof  is the discovery of a 
new monotone quantity associated to the deformed 
metric. We divide the proof into two subsections.

\subsection{A monotonicity property for quasi-spherical 
metrics on a Schwarzschild background} \label{gmono}

Consider part of a spatial Schwarzschild manifold 
\be
( M^m_{r_0} , g^m ) = 
\left( S^2 \times [ r_0 , \infty) , \frac{1}{ 1 - \frac{ 2m }{ r } } 
dr^2 + r^2 d \sigma^2 \right),
\ee
where $ r_0 $ is a constant chosen to satisfy 
$ \left\{ 
\begin{array}{ll}
r_0  >  2m, & \mathrm{if } \ m \geq 0 \\
r_0  >  0, & \mathrm{if}  \ m < 0 .
\end{array}
\right.
$ 
Here $ m $ is the ADM mass of the Schwarzschild 
metric $  g^m $, $ r $ is the radial coordinate on 
$ [r_0, \infty) $, and $ d \sigma^2 $ denotes the 
standard metric on the unit sphere $ S^2 \subset \R^3 $. 

Let $ N $ be the positive function on $ M^m_{r_0} $ 
defined by
\be \label{dfofn}
N = \sqrt{ 1 - \frac{ 2m }{ r } } . 
\ee
In terms of $ N $, $ g^m $ takes the form
\be
g^m = \frac{1}{N^2} dr^2 + r^2 d \sigma^2 .
\ee
The next lemma follows directly from the existence 
theory established in \cite{Shi-Tam} 
(see also \cite{Bartnik}). 

\begin{lma} \label{schwarzext}
Let $ \S_0 $ be the boundary of $ (M^m_{r_0}, g^m)$.
Given any positive function $ \phi $ on $ \S_0$,
there exists a positive function $ u $ on $ \Mmrz $ 
such that 
\begin{enumerate}
\item[(i)] 
The metric 
\be
g^m_u = \left( \frac{ u }{ N } \right)^2 d r^2 + r^2 d \sigma^2 
\ee
has zero scalar curvature and is asymptotically flat . 
\item[(ii)] The mean curvature of $ \S_0 $ (with respect to the 
$ \infty $-pointing normal) in $ (\Mmrz, g^m_u ) $ is equal to $ \phi $.
\item[(iii)] The quotient
$ \frac{ u }{ N } $ has the asymptotic expansion
\be
 \frac{ u }{ N }  = 1 + \frac{ m_0 }{ r } 
 + O\left( \frac{1}{r^2} \right)  \ 
 \mathrm{as} \ r \rightarrow \infty ,
\ee
where $ m_0 $ is the ADM mass of $ g^m_u $. 
\end{enumerate} 

\end{lma}

\begin{proof} Consider a Euclidean background metric
\be
ds^2 = dr^2 + r^2 d \sigma^2 
\ee
on $ \Mmrz = S^2 \times [r_0, \infty) $. By Theorem 2.1 
in \cite{Shi-Tam}, there is a unique
positive function $ v $ on $ \Mmrz$ such that 
\be
g_v = v^2 d r^2 + r^2 d \sigma^2
\ee
has zero scalar curvature, is asymptotically flat and
the mean curvature of $ \S_0 $ in $ (\Mmrz, g_v) $ 
is given by $ \phi $. 
Furthermore, $ v $ has
an asymptotic expansion
\be \label{asympofv}
v = 1 + \frac{ m_0 }{ r } + O\left( \frac{1}{r^2} \right) ,
\ee
where $ m_0 $ is the ADM mass of $ g_v $.
Let $ u = N v $, Lemma \ref{schwarzext} is proved.  \end{proof}

We note that metrics of the form
$ v^2 dr^2 + r^2 d \sigma^2 $
are called (shear free) {\em quasi-spherical} 
metrics \cite{Bartnik}.
By the formula (2.26) in \cite{Bartnik} (or 
 (1.10) in \cite{Shi-Tam}), 
the differential equation satisfied by $ v = u N^{-1} $
in Lemma {\ref{schwarzext} is 
\be \label{ueq}
\frac{ 2 }{ r }  \frac{ \partial v }{ \partial r }  = 
\frac{ v^2 }{ r^2 } \triangle_{ S^2 } v +
 \frac{( v- v^3)}{ r^2 },
\ee
where $ \triangle_{ S^2 } $ denotes the Laplacian 
operator of  the metric $d \sigma^2 $ on $ S^2 $.

\begin{prop} \label{newmonotone}
Let $ u $, $ g^m_u $, $ m_0 $ be given in Lemma \ref{schwarzext}. 
Let $ \Sr $ be the radial coordinate sphere in 
$ \Mmrz $, i.e. $ \Sr = S^2 \times \{ r \} $. 
Let $ H_{_S} $, $ H_u $ be the mean curvature 
of $ \Sr $ with respect to the metric
$ g^m $, $ g^m_u$. Then
$$ \oint_{ \S_r } N ( H_{_S} - H_u )  \ d \sigma $$ 
is monotone non-increasing in $ r $. 
Furthermore,
\be \label{limit}
\lim_{ r \rightarrow \infty } \oint_{ \Sr } N ( H_{_S} - H_u )  \
 d \sigma  = 8 \pi ( m_0 - m ) . 
\ee
\end{prop}

\begin{proof} We have
$
H_{_S} = \frac{ 2 }{r} N$ {and} $  H_u =   
\frac{ 2 }{r} v^{-1} ,
$
where $ v = u N^{-1} $.  
Hence 
\beq
 \oint_{ \Sr } N ( H_{_S} - H_u )  \ d \sigma  & = & 
  \oint_{ \Sr } \left( \frac{ 2 }{ r } \right) 
 ( N^2 - N v^{-1} )  \ d \sigma  \nonumber \\
 & = & 
 \oint_{ S^2 }  2 r ( N^2 - N v^{-1} )  \ d \omega, 
\eeq
where $ d \omega = r^{-2} d \sigma $ is the surface 
measure  of $ d \sigma^2 $ on $ S^2 $. 
As $
N^2 = 1 - \frac{ 2m }{ r },
$
we have
\beq \label{simpfy}
 \oint_{ \Sr } N ( H_{_S} - H_u )  \ d \sigma  & = & 
 \oint_{ S^2 }  ( 2 r   -  4 m  - 2 r N v^{-1} )  \ d \omega . 
\eeq
Therefore,
\beq \label{drho}
 \frac{d}{d r } \oint_{\Sr } N  ( H_{_S} - H_u )  \ d \sigma 
& =  & \oint_{ S^2 }  \left[ ( 2  - 2 N v^{-1} )  
 -  2 r \frac{\partial N}{\partial r} v^{-1} \right] \ d \omega \nonumber \\
& & +  \oint_{S^2}
2 r N v^{-2} \frac{\partial v}{\partial r}  \ d \omega .
\eeq
By (\ref{ueq}), we have
\be 
 v^{-2} \frac{ \partial v }{ \partial r }  = 
 \frac{ 1}{ 2r } \triangle_{ S^2 } v +
 \frac{( v^{-1}- v)}{ 2 r } . 
\ee
Thus the last term in (\ref{drho}) becomes
\beq
\oint_{ S^2 } 2 r N v^{-2} \frac{\partial v}{\partial r}  
\ d \omega   & = & \oint_{ S^2 }   N \triangle_{ S^2 } v 
\ d \omega + \oint_{ S^2 } N ( v^{-1}  - v )  \ d \omega 
\nonumber \\
& = &  \oint_{ S^2 } N ( v^{-1}  - v ) \  d \omega ,
\eeq
where we have used the fact that $ N $ is 
a constant on  each $ \Sr $
and $ \oint_{S^2} \triangle_{S^2} v \ d \omega  = 0 $.
Hence the right side of (\ref{drho}) is given by 
\beq \label{aa}
 \oint_{ S^2 }  \left[ ( 2  - 2 N v^{-1} )  
 -  2 r \frac{\partial N}{\partial r} v^{-1} 
 +  N ( v^{-1}  - v )  \right]  \ d \omega  . 
 \eeq
 Replace $ v $ by $ u N^{-1} $, the integrand of 
 (\ref{aa}) becomes
 \be \label{integrand}
  2  -  N^2 u^{-1} 
 -  2 r \frac{\partial N}{\partial r} N u^{-1} - u .
 \ee
By \eqref{dfofn}, 
we have
 \be \label{algebra}
 N^2 + 2 r N \frac{\partial N}{\partial r}  = 1 . 
 \ee
 Therefore, it follows from (\ref{drho}), (\ref{aa}), 
 (\ref{integrand}) and  (\ref{algebra}) that
\be \label{derivativeformula}
\frac{d}{d r } \oint_{\Sr } N  ( H_{_S} - H_u )  \ d \sigma 
= - \oint_{ S^2} u^{-1} ( u - 1)^2 \ d \omega ,
\ee
which proves that $ \oint_{\Sr} N  ( H_{_S} - H_u ) \  d \sigma $
is monotone non-increasing in $ r $.

To evaluate $
\lim_{ r \rightarrow \infty} \oint_{\Sr } N ( H_{_S} - H_u) \ d \sigma ,
$
we have 
\be
N v^{-1} = 1 - \frac{ ( m_0 + m ) }{ r } 
+ O \left( \frac{ 1 }{ r^2 } \right) 
\ee
by (\ref{dfofn}) and (\ref{asympofv}).
Therefore, by (\ref{simpfy}) we have
\be
\oint_{ \Sr } N ( H_{_S} - H_u )  \ d \sigma   =
 \oint_{ S^2 }   2 ( m_0   -  m )  \ d \omega 
 +  O ( r^{-1} ) ,
\ee
which implies 
\be
\lim_{ r \rightarrow \infty} \oint_{\Sr } N ( H_{_S} - H_u) 
 d \sigma  = 8 \pi ( m_0 - m) . 
\ee
Proposition \ref{newmonotone} is proved. \end{proof}

\subsection{Application of the Riemannian Penrose 
Inequality} \label{apofrpi}

In this section, we  glue a body surrounding horizons, 
whose outer boundary is metrically a round sphere, 
to an asymptotically flat manifold $ (\Mmrz, g^m_u )$ 
constructed in Lemma \ref{schwarzext}, and apply 
the Riemannian Penrose Inequality and 
Proposition \ref{newmonotone} 
to prove Theorem \ref{mthm}. 

We start with the following lemma.

\begin{lma} \label{gmtlma}
Let $ \Om $ be a body surrounding horizons. 
Suppose its outer boundary $ \SO $ 
has positive mean curvature, then its horizon 
boundary $ \SH $ strictly 
minimizes area among all closed surfaces in 
$ \Om $ that enclose $ \SH $.
\end{lma}

\begin{proof}  As $ \Om $ is compact and the mean curvature 
vector of $ \SO $ points into $ \Om$, 
it follows from the standard geometric measure theory
that there exist surfaces that minimize area among 
all closed surfaces  in $ \Om $
that enclose $ \SH $, furthermore none of the 
minimizers touches $ \SO $. 
Let $ \S $ be any such a minimizer. 
By the Regularity Theorem 1.3  in \cite{IMF}, 
$ \S $ is a $ C^{1,1} $ surface, and is $ C^\infty $ 
where it does not touch $ \SH $; moreover, the mean 
curvature of $ \S $ is $ 0 $ on 
$ \S \setminus \SH $ and equals the 
mean curvature  of $ \SH $ 
$ \mathcal{H}^2$-a.e. on $ \S \cap \SH$. 
Suppose $ \S $ is not identically $ \SH $. 
As $ \SH $ has zero mean curvature,  the
maximum principle implies that $ \S $ does not touch $ \SH $. 
Hence, $ \S $ is a smooth closed minimal surface in the 
interior of $ \Om $,  contradicting the assumption
that $ \Om $ has no other closed minimal surfaces except 
$ \SH $. Therefore, $ \S $ must be identically $ \SH $. \end{proof}


Let $ \Om $ be a body surrounding horizons given in 
Theorem \ref{mthm}. Let $ R $ and $ R_H $ be the area 
radii of $ \SO $ and $ \SH $, which are defined by
\be
4 \pi R^2 = | \SO | \ \ \mathrm{and} \ \ 4 \pi R_H^2 = | \SH | . 
\ee
It follows from  Lemma \ref{gmtlma} that 
$
R > R_H .
$
To proceed, we choose $ ( M^m , g^m ) $ to be one-half 
of a spatial Schwarzschild manifold whose horizon  
has the same area as $ \SH $, i.e.
\be
( M^m , g^m ) = \left( S^2 \times [ R_H , \infty) , 
\frac{1}{ 1 - \frac{ 2m }{ r } } dr^2 + r^2 
d \sigma^2 \right),
\ee
where $ m $ is chosen to satisfy $ 2 m = R_H $. 
As $ R >  R_H$, $ \SO $ can be isometrically 
imbedded in $ \SMm $ as the coordinate sphere 
\be
\S_R = \{  r = R \} . 
\ee
Henceforth, we identify $ \SO $ with $ \S_R $ 
through this isometric imbedding.
Let $ \MmR $ denote the exterior of $ \SR $ in $ M^m $.
By Lemma \ref{schwarzext} and 
Proposition \ref{newmonotone},  there exists a metric 
\be
g^m_u = \left( \frac{ u }{ N } \right)^2 d r^2 + r^2 d \sigma^2 
\ee
on $ \MmR $ such that $ g^m_u $ has zero scalar 
curvature, is asymptotically flat, and the mean curvature 
of $ \SR $  (with respect to the $ \infty $-pointing normal) 
in $ (\MmR, g^m_u) $ agrees with $ H $, the mean 
curvature of $ \SO $  in $ \Om$.
Furthermore, the integral 
\be
\oint_{\Sr} N ( H_{_S} - H_u ) \ d \sigma 
\ee
is monotone non-increasing in $ r $ and converges to 
$ 8 \pi ( m_0 - m ) $ as $ r \rightarrow \infty $, 
where $ m_0 $ is the ADM mass of $ g^m_u $. 

Now we attach this asymptotically flat manifold 
$ ( \MmR, g^m_u )$  to the compact body $ \Om $ 
along $ \SO $ to get a complete Riemannian manifold 
$ M $ whose boundary is $ \SH $. The resulting metric 
$ \gm $ on $ M $ satisfies the properties that it is 
asymptotically flat, has nonnegative scalar curvature 
away from $ \SO $,  is Lipschitz near $ \SO $, and the 
mean curvatures of $ \SO $ computed in both sides 
of $ \SO $ in $ M $  (with respect to the $\infty$-pointing normal)  
agree identically. 

\begin{lma} \label{bdryminimizing} 
The horizon boundary $ \SH $ is strictly outer minimizing 
in $ M $, i.e $ \SH $ strictly minimizes area among all 
closed surfaces in $ M $ that enclose $ \SH $.
\end{lma}

\begin{proof} By the construction of $ g^m_u $,  we know $ (\MmR, g^m_u ) $ 
is foliated by  $\{ \Sr \}_{ r \geq R} $, where each $ \Sr $ 
has positive mean curvature. Let $ \S $ be a surface 
that minimizes area among surfaces in $ M $ that 
encloses $ \SH $ (such a minimizer exists as $ M $ 
is asymptotically flat). We claim that  $ \S \setminus \Om $ 
must be empty, for otherwise  $ \S \setminus \Om $
would be a smooth, compact  minimal surface in 
$ (\MmR, g^m_u) $ with boundary lying in $ \SO $,  
and that would contradict the maximum principle. 
Therefore, $ \S \subset \Om $. It then follows 
from Lemma \ref{gmtlma} that $ \S  = \SH $.  \end{proof}


The next lemma is an application of the ``corner smoothing" 
technique  in \cite{Miao_PMT}.

\begin{lma} \label{reflectionandsmoothing}
There exists a sequence of smooth asymptotically flat metrics 
$ \{ h_k \} $ defined  on the background manifold of $ M $ 
such that $ \{ h_k \} $ converges uniformly to $ \gm $ 
in the $ C^0 $ topology, each $ h_k $ has nonnegative 
scalar curvature, $ \SH $ has zero mean curvature 
with respect to each $ h_k $ (in fact $ \SH $ can be 
made totally geodesic w.r.t $ h_k $),  
and the ADM mass of $ h_k $ converges to the 
ADM mass of $ \gm $.
\end{lma}

\begin{proof}   
Let $ M^\prime $ be an exact copy of $ M $. We glue 
$ M $ and $ M^\prime $ along their common boundary 
$ \SH $ to get a Riemannian 
manifold $ \bN $ with two asymptotic ends. 
Let $ g_{_{\bN}} $ be the resulting metric on $ \bN $ and let 
$ \SO^\prime $  be the copy of $ \SO $ in $ M^\prime $. 
Denote by $ \S $ the union of  $ \SO $, $ \SH $ and 
$ \SO^\prime $, we then know that the mean curvatures 
of $ \S $ computed in both sides of $ \S $ in $ \bN $ 
(with respect to normal vectors pointing to the same end of 
$ \bN $) agree. (At $ \SO $ and $ \SO^\prime$, this is 
guaranteed by the construction of $ g^m_u $, and 
at $ \SH $, this is provided by the fact that $ \SH $ 
has zero mean curvature.)

Apply Proposition 3.1 in \cite{Miao_PMT} to $ \bN $ 
at $ \S $, followed by a conformal deformation as 
described in Section 4.1 in \cite{Miao_PMT}, we get 
a sequence of smooth asymptotically flat metrics $ \{ g_k \} $, 
defined on the background manifold of $ \bN $, with 
nonnegative scalar curvature such that  $ \{ g_k \} $ 
converges uniformly to $ g_{_{\bN}} $ in the $ C^0 $ 
topology and the ADM mass of $ g_k $ converges to the 
ADM mass of $ g_{_{\bN}} $ on both ends of $ \bN $. 
Furthermore, as $ \bN $ has a reflection isometry 
(which maps a point $ x \in M $ to its copy in $ M^\prime $), 
detailed checking of the construction in Section 3 
in \cite{Miao_PMT} shows that $ \{ g_k \} $ can be 
produced in such a way that each $ g_k $ also has the same 
reflection isometry. (Precisely, this can be achieved by 
choosing the mollifier $ \phi(t) $ in equation (8) 
in \cite{Miao_PMT} and the cut-off function $ \sigma(t) $ 
in equation (9) in \cite{Miao_PMT} to be both even functions.)
Therefore, if we let $ \bN_k $ be the Riemannian manifold 
obtained by replacing the metric $ g_{_{\bN}} $ by $ g_k $
on $ \bN$, then $ \SH $ remains a surface with zero mean 
curvature in $ \bN_k $ (in fact $\SH$ is  totally geodesic).  
Define $ h_k $ to be the restriction of $ g_k $ to the 
background manifold of $ M $, Lemma \ref{reflectionandsmoothing} 
is proved. \end{proof}


We continue with the proof of Theorem \ref{mthm}.
Let $ \{ h_k \} $ be the metric approximation of $ \gm $ 
provided in Lemma \ref{reflectionandsmoothing}.
Let $ M_k $ be the asymptotically flat manifold obtained 
by replacing the metric $ \gm $ on $ M $ by $ h_k $. 
For any surface $ \tilde{\Sigma} $ in $ M $, let 
$ | \tilde{\Sigma} |_k $, $ | \tilde{\Sigma} |$
 be the area of $ \tilde{\Sigma} $ 
w.r.t the induced metric from $ h_k $, $ \gm $ respectively. 
We can not apply the Riemannian Penrose Inequality 
directly to claim 
$ m_{ADM}(h_k ) \geq \sqrt{ \frac{ | \SH |_k }{ 16 \pi } } $. 
That is because we do not know if $ \SH $ remains
to be the outermost minimal surface in $ M_k $. However, 
since $ \SH $ is a minimal surface in $ M_k $, 
we know the outermost minimal surface  in $ M_k $, 
denoted by $ \S_k $, exists and its area satisfies 
\be \label{areask}
 | \Sigma_k |_k = \inf \{ | \tilde{\Sigma} |_k  \ | \ 
 \tilde{\Sigma} \in \mathcal{S} \}
\ee
where $ \mathcal{S} $ is the set of  closed surfaces 
$ \tilde{\Sigma} $
 in $ M $ that enclose $ \SH $ (see \cite{Bray}, \cite{IMF}).
By the Riemannian Penrose  Inequality (Theorem 1 in \cite{Bray}), 
we have
\be \label{approxrpi}
m_{ADM}(h_k ) \geq \sqrt{ \frac{ | \Sigma_k |_k }{ 16 \pi } }.
\ee
Let $ k $ approach infinity, we have
\be \label{approxmass}
\lim_{ k \rightarrow \infty} m_{ADM} (h_k) = m_{ADM} (\gm) , 
\ee
and
\be \label{approxarea}
\lim_{ k \rightarrow \infty} | \Sigma_k |_k  = 
\inf \{ | \tilde{\Sigma} |  \ | \ 
 \tilde{\Sigma} \in \mathcal{S} \} 
\ee
where we have used \eqref{areask} and the fact that $ \{ h_k \} $ converges uniformly to $ g_{_M} $ in the 
$ C^0 $ topology. 
By Lemma \ref{bdryminimizing},  we also have
\be \label{areaofsh}
| \Sigma_H | = \inf \{ | \tilde{\Sigma} |  \ | \ 
 \tilde{\Sigma} \in \mathcal{S} \} .
\ee 
Therefore,  it follows from (\ref{approxrpi}), (\ref{approxmass}),
(\ref{approxarea})  and \eqref{areaofsh}  that
\be \label{nonsmoothrpi}
m_{ADM} (\gm) \geq \sqrt{ \frac{ | \SH | }{ 16 \pi } } . 
\ee

To finish the proof of Theorem \ref{mthm}, we make use 
of the  monotonicity of the integral 
\be
\oint_{\S_r} N ( H_{_S} - H_u ) \ d \sigma . 
\ee
By Proposition \ref{newmonotone}, we have
\beq \label{appofmonotone}
\oint_{\SO} N ( H_{_S} - H_u ) \ d \sigma & \geq & 
\lim_{ r \rightarrow \infty} \oint_{\S_r} N ( H_{_S} - H_u ) 
\ d \sigma \nonumber \\
& = & 8 \pi ( m_0 - m) .
\eeq
On the other hand, we know 
\be \label{mzmassn}
m_0 =  m_{ADM}(g^m_u ) = m_{ADM} (g_M) ,
\ee
and
\be \label{hchm}
m = \frac{ 1}{2} R_H  = \sqrt{ \frac{ | \SH | }{ 16 \pi } } .
\ee
Therefore, it follows from (\ref{appofmonotone}), 
(\ref{mzmassn}), (\ref{hchm}) and (\ref{nonsmoothrpi}) that 
\be \label{bdryineq}
\oint_{\SO} N ( H_{_S} - H_u ) d \sigma \geq 0 . 
\ee
Plug in 
$H_{_S} = \frac{ 2 }{ R } N$, $ H_u = H$
and $ N  = \sqrt{ 1 - \frac{ R_H }{ R } } , $
we then have
\be \label{bridge}
8 \pi  R \sqrt{ 1 - \frac{ R_H }{ R } } \geq \oint_{\SO} H \ d \sigma .
\ee
Direct computation shows that (\ref{bridge}) 
is equivalent to (\ref{mineq}). Hence, (\ref{mineq}) is proved.

Finally, when the equality in (\ref{mineq}) holds, we have
\be
\oint_{\Sr} N ( H_{_S} - H_u ) d \sigma  = 0 , \ \ \forall \ r \geq R . 
\ee
By the derivative formula (\ref{derivativeformula}), 
 $ u $  is identically $ 1 $ on $ \MmR $. Therefore, 
the metric $ g^m_u $ is indeed the Schwarzschild metric 
$ g^m $. Since the mean curvature of $ \SO $ in 
$ (\MmR, g^m_u ) $ was arranged to equal $H $, the mean 
curvature of $ \SO $ in $ \Om$,  we conclude that  
$ H = \frac{ 2 }{ R} \left( 1 - \frac{ R_H}{R} \right)^{\frac{1}{2}}  $,
which is a constant. 
Theorem \ref{mthm} is proved. 
\stop

\vh

Comparing to the equality case in Theorem 2,  one would expect that the equality in  \eqref{mineq} holds if and only if  $ \Om $ is  isometric to a region, in a spatial Schwarzschild manifold, which is bounded by a rotationally symmetric sphere and the Schwarzschild horizon. We believe that this is true, but are not able  to prove it at this stage.  A confirmation of this expectation  seems to require a good knowledge of the behavior of a  sequence of asymptotically flat $ 3 $-manifolds with  controlled $ C^0$-geometry,  on which the equality of  the Riemannian Penrose Inequality is nearly satisfied.  We leave this as an open question.

\section{Some discussion} \label{rmk} 
Let $ \Sigma $ be an arbitrary closed $ 2 $-surface in a general $ 3 $-manifold $ M $ with nonnegative scalar curvature.  
Consider the quantity 
\be  \label{dfofm}
m (\Sigma) =  \sqrt{ \frac{ | \Sigma | }{ 16 \pi } } 
\left[ 1 - \frac{ 1 }{ 1 6 \pi | \Sigma | } \left( \oint_{\Sigma} H 
\right)^2 \right] 
\ee
where $ | \Sigma | $ is the area of $ \Sigma $, $ H $ is the 
mean curvature of $ \Sigma $ in $ M $ and we omit the 
surface measure $ d \sigma $ in the integral.
Theorem \ref{mthm} suggests that,
if $ \Sigma $ is metrically a round sphere, 
 $ m ( \Sigma ) $ may potentially agree with a hidden definition of quasi-local mass  of  $ \Sigma $.
Such a speculation could be further strengthened by the resemblance between 
$ m (\Sigma )$ and the  Hawking quasi-local mass \cite{Hawking}
\be m_H (\Sigma) =  \sqrt{ \frac{ | \Sigma | }{ 16 \pi } } 
\left[ 1 - \frac{ 1 }{ 1 6 \pi  }  \oint_{\Sigma} H^2   \right].
\ee
By H\"{o}lder's inequality, we have
\be \label{mandmh}
 m ( \Sigma ) \geq m_H ( \Sigma) 
\ee
for any surface $ \Sigma $. 
On the other hand, if $ \S $ is a closed convex surface in the Euclidean space
$ \R^3 $, the classic Minkowski inequality \cite{Polya_Szego}
\be \label{Mkineq}
\left( \oint_\Sigma H \ d \sigma \right)^2 \geq 16 \pi | \Sigma |
\ee
 implies that  $ m (\Sigma ) \leq 0 $ and $ m ( \Sigma ) = 0 $ 
if and only if 
$ \Sigma $ is a round sphere in $ \R^3 $. Therefore, even 
though bigger than $ m_H ( \Sigma ) $, 
$ m ( \Sigma ) $ shares the same character as 
$ m_H ( \Sigma ) $ that it is  negative on most 
convex surfaces in $ \R^3 $. 

In order to gain positivity and to maintain the same numerical 
value on metrically round spheres, we propose to modify $ m ( \Sigma ) $ in a similar way as the Brown-York mass $ m_{BY} ( \Sigma ) $ \cite{BYmass1} is defined. 
Recall that, for those $ \Sigma $ with positive Gaussian
curvature,  $ m_{BY} ( \Sigma ) $ is defined to be
\be
m_{BY} ( \Sigma ) = \frac{ 1 }{ 8 \pi } 
\left( \oint_{ \Sigma } H_0 \ d \sigma - \oint_{ \Sigma } H \ d \sigma \right) 
\ee
where $ H_0 $ is the  mean curvature of $ \Sigma $ when 
it is isometrically embedded in $ \R^3 $. 
Now suppose $ \Sigma $ is  metrically a round sphere, 
then
\be
\left( \oint_\Sigma H_0 \right)^2 = 16 \pi | \Sigma |.
\ee
In this case, we can re-wriite $ m ( \Sigma ) $  as either
\be
m ( \Sigma ) =   \sqrt{ \frac{ | \Sigma | }{ 16 \pi } } 
\left[ 1 -  \left(  \frac{ \oint_{\Sigma} H }{ \oint_\Sigma H_0  } \right)^2 \right] 
\ee
or 
\be
m ( \Sigma ) =    \frac{ 1 }{ 16 \pi } \left( \oint_\Sigma H_0 \right)
\left[ 1 -  \left(  \frac{ \oint_{\Sigma} H  }{ \oint_\Sigma H_0  } \right)^2 \right] .
\ee
This motivates us to consider the following two quantities: 

 \begin{df}
 For any $ \Sigma $ with positive Gaussican curvature, define
\be
m_1 ( \Sigma ) = \sqrt{ \frac{ | \Sigma | }{ 16 \pi } } 
\left[ 1 -  \left(  \frac{ \oint_{\Sigma} H  }{ \oint_\Sigma H_0 } \right)^2 \right], 
\ee
and
\be
m_2 ( \Sigma) = \frac{ 1 }{ 16 \pi } \left( \oint_\Sigma H_0 \right)
\left[ 1 -  \left(  \frac{ \oint_{\Sigma} H  }{ \oint_\Sigma H_0 } \right)^2 \right] ,
\ee
 where $ H $ is the mean curvature of $ \Sigma $ in $ M $ 
 and  $ H_0 $ is the  mean curvature of $ \Sigma $ when it is isometrically embedded in $ \R^3 $.
\end{df}

The following result compares $ m_H ( \Sigma ) $, 
$ m_1 ( \Sigma) $, $ m_2 ( \Sigma ) $ and
$ m_{BY} ( \Sigma )$. 

\begin{thm} \label{qlmassthm}
Suppose $ \Sigma $ is a closed $ 2 $-surface with positive Gaussian
curvature in a $ 3 $-manifold $ M $. Then 
\begin{enumerate}
\item[(i)] $ m_1 ( \Sigma ) \geq m_H ( \Sigma ) $, and equality holds if and only if $ \Sigma $ is  metrically a round sphere and $ \Sigma $ has constant mean curvature. 

\item[(ii)]  $ m_{BY} ( \Sigma )  \geq  m_2 ( \Sigma ) $, and  equality holds 
if and only if $ \oint_\Sigma H_0 \ d \sigma = \oint_\Sigma H \ d \sigma $.

\item[(iii)] Suppose $ \Sigma $ bounds a domain $ \Omega $ with nonnegative scalar  curvature and the mean curvature of $\Sigma $ 
in $ \Omega $ is positive, then 
$$  m_2 ( \Sigma ) \geq m_1 ( \Sigma ) \geq 0 .$$
Moreover,  $ m_1 ( \Sigma ) = 0 $ if and only if $ \Omega $ is isometric to a domain in $ \R^3 $, 
and $ m_2 ( \Sigma ) = m_1 ( \Sigma ) $ if and only if either 
$ \Omega $ is isometric to a domain in $ \R^3 $ in which case
$ m_2 ( \Sigma ) = m_1 ( \Sigma ) = 0 $ or 
$ \Sigma $ is metrically a round sphere.

\end{enumerate}
\end{thm}

\begin{proof} (i) Let $ m (\Sigma ) $ be defined as in \eqref{dfofm}.  
By the Minkowski inequality \eqref{Mkineq}, we have
$ m_1 ( \Sigma ) \geq m ( \Sigma ).  $
By  \eqref{mandmh}, we have
$ m ( \Sigma ) \geq m_H ( \Sigma).$
Therefore,
$ m_1 ( \Sigma ) \geq m_H ( \Sigma ) $
and equality holds if and only if $ \Sigma $ is metrically a round sphere and the mean curvature  of $ \Sigma $ in $ M $ is a constant. 

(ii) This case is elementary. Let $ a = \oint_\Sigma H $ and 
$ b = \oint_\Sigma H_0 $. Then (ii) is equivalent to the 
inequality $ \left( 1 - \frac{a}{b} \right)^2 \geq 0 $. 

(iii) By the result of Shi and Tam \cite{Shi-Tam}, i.e. Theorem \ref{stthm}, we have
\be \label{rest}
1 -  \left(  \frac{ \oint_{\Sigma} H  }{ \oint_\Sigma H_0 } \right)^2 
\geq 0 
\ee
with equality holding if and only if $ \Omega $ is isometric to 
a domain in $ \R^3 $. (iii) now follows directly from \eqref{rest}
and the Minkowski inequality \eqref{Mkineq}. 
\end{proof}

Suppose $ \Omega $ is a compact $ 3 $-manifold with boundary 
with nonnegative scalar curvature and its boundary $ \partial \Omega $ has positive Gaussian curvature and positive mean curvature.
Theorem \ref{qlmassthm} implies that 
\be
m_{BY} ( \partial \Omega) \geq m_2 ( \partial \Omega) 
\geq m_1 ( \partial \Omega ) \geq m_H ( \partial \Omega )
\ee
with $ m_1 ( \partial \Omega) \geq 0 $ 
and $ m_{BY} ( \partial \Omega ) = m_H ( \partial \Omega ) $
if and only if $ \Omega $ is isometric to a round ball in $ \R^3 $. 
This provides a slight generalization of a previous result of Shi and Tam (Theorem 3.1 (b) in \cite{Shi-Tam07}), which showed 
$ m_{BY} ( \partial \Omega ) \geq m_H( \partial \Omega)$.

\bibliographystyle{plain}

\begin{thebibliography}{10}

\bibitem{ADM}
R.~Arnowitt, S.~Deser, and C.~W. Misner.
\newblock Coordinate invariance and energy expressions in general relativity.
\newblock {\em Phys. Rev. (2)}, 122:997--1006, 1961.

\bibitem{Bartnik_qlmass}
Robert Bartnik.
\newblock New definition of quasilocal mass.
\newblock {\em Phys. Rev. Lett.}, 62(20):2346--2348, 1989.

\bibitem{Bartnik}
Robert Bartnik.
\newblock Quasi-spherical metrics and prescribed scalar curvature.
\newblock {\em J. Differential Geom.}, 37(1):31--71, 1993.


\bibitem{Bray_Thesis}
Hubert~L. Bray.
\newblock The {P}enrose inequality in general relativity and volume comparison
  theorems involving scalar curvature.
\newblock {\em Stanford University Thesis}, 1997.

\bibitem{Bray}
Hubert~L. Bray.
\newblock Proof of the {R}iemannian {P}enrose inequality using the positive
  mass theorem.
\newblock {\em J. Differential Geom.}, 59(2):177--267, 2001.

\bibitem{BYmass1}
J.~David Brown and James~W. York, Jr.
\newblock Quasilocal energy in general relativity.
\newblock In {\em Mathematical aspects of classical field theory (Seattle, WA,
  1991)}, volume 132 of {\em Contemp. Math.}, pages 129--142. Amer. Math. Soc.,
  Providence, RI, 1992.


\bibitem{CDO-Yau}
D.~Christodoulou and S.-T. Yau.
\newblock Some remarks on the quasi-local mass.
\newblock In {\em Mathematics and general relativity (Santa Cruz, CA, 1986)},
  volume~71 of {\em Contemp. Math.}, pages 9--14. Amer. Math. Soc., Providence,
  RI, 1988.



\bibitem{Hawking}
Stephen Hawking.
\newblock Gravitational radiation in an expanding universe.
\newblock {\em J. Mathematical Phys.}, 9:598--604, 1968.

\bibitem{IMF}
Gerhard Huisken and Tom Ilmanen.
\newblock The inverse mean curvature flow and the {R}iemannian {P}enrose
  inequality.
\newblock {\em J. Differential Geom.}, 59(3):353--437, 2001.


\bibitem{Miao_PMT}
Pengzi Miao.
\newblock Positive mass theorem on manifolds admitting corners along a
  hypersurface.
\newblock {\em Adv. Theor. Math. Phys.}, 6(6):1163--1182 (2003), 2002.

\bibitem{Nirenberg}
Louis Nirenberg.
\newblock The {W}eyl and {M}inkowski problems in differential geometry in the
  large.
\newblock {\em Comm. Pure Appl. Math.}, 6:337--394, 1953.


\bibitem{Polya_Szego}
G.~P{\'o}lya and G.~Szeg{\"o}.
\newblock {\em Isoperimetric {I}nequalities in {M}athematical {P}hysics}.
\newblock Annals of Mathematics Studies, no. 27. Princeton University Press,
  Princeton, N. J., 1951.


\bibitem{Schoen-Yau}
Richard Schoen and Shing~Tung Yau.
\newblock On the proof of the positive mass conjecture in general relativity.
\newblock {\em Comm. Math. Phys.}, 65(1):45--76, 1979.

\bibitem{Shi-Tam}
Yuguang Shi and Luen-Fai Tam.
\newblock Positive mass theorem and the boundary behaviors of compact manifolds with nonnegative scalar curvature.
\newblock {\em J. Differential Geom.}, 62(1):79--125, 2002.

\bibitem{Shi-Tam07}
Yuguang Shi and Luen-Fai Tam.
\newblock Quasi-local mass and the existence of horizons.
\newblock {\em Comm. Math. Phys.}, 274(2):277--295, 2007.


\bibitem{Witten}
Edward Witten.
\newblock A new proof of the positive energy theorem.
\newblock {\em Comm. Math. Phys.}, 80(3):381--402, 1981.


\end{thebibliography}

\end{document}